\documentclass[12pt]{article}

\usepackage{mathrsfs} 

\usepackage[all]{xy}

\newtheorem{tm}{Theorem}[subsection]
\newtheorem{lm}[tm]{Lemma}
\newtheorem{pr}[tm]{Proposition}
\newtheorem{rmk}[tm]{Remark}

\newtheorem{ex}[tm]{Example}

\newtheorem{??}[tm]{Question}
\newtheorem{defi}[tm]{Definition}

\newtheorem{ass}[tm]{Assumption}

\font\tenmsb=msbm10
\font\sevenmsb=msbm7
\font\fivemsb=msbm5
    
\newfam\msbfam
\textfont\msbfam=\tenmsb
\scriptfont\msbfam=\sevenmsb
\scriptscriptfont\msbfam=\fivemsb
\def\Bbb#1{{\fam\msbfam #1}}

\font\teneufm=eufm10
\font\seveneufm=eufm7
\font\fiveeufm=eufm5
\newfam\eufmfam
\textfont\eufmfam=\teneufm
\scriptfont\eufmfam=\seveneufm
\scriptscriptfont\eufmfam=\fiveeufm
\def\frak#1{{\fam\eufmfam\relax#1}}

\def\lorw{\longrightarrow}
\newcommand\n{\noindent}

\newcommand\ci{\cite}
\newcommand\s{\sigma}
\newcommand\rat{{\Bbb Q}}
\newcommand\comp{{\Bbb C}}

\newcommand\zed{{\Bbb Z}}

\newcommand\pn[1]{{\Bbb P}^{#1}}

\newcommand\blacksquare{{\hspace*{\fill} $\fbox{}$}}

\newcommand\e{\epsilon}

\newcommand{\phix}[2]{ \,^p\!{\cal H}^{#1}(#2)}

\newcommand{\im}{ \hbox{\rm Im} }
\newcommand{\coim}{ \hbox{\rm Coim} }

\newcommand{\ke}{ \hbox{\rm Ker} }

\newcommand{\ben}{\begin{enumerate}}
\newcommand{\een}{\end{enumerate}}

\newcommand{\bit}{\begin{itemize}}
\newcommand{\eit}{\end{itemize}}

\newcommand{\beq}{\begin{equation}}
\newcommand{\eeq}{\end{equation}}

\newcommand{\la}{\label}

\title{Hodge-theoretic splitting mechanisms\\for projective  maps\\in appendix: a letter  from P. Deligne}
\author{
Mark Andrea A.  de Cataldo
}
\date{2013}

\begin{document}\maketitle

\begin{abstract}
According to the decomposition and relative hard Lefschetz  theorems, given a projective map
of complex quasi projective  algebraic varieties and a relatively ample line bundle, the rational  intersection cohomology groups of the domain of the map
split  into various direct summands. While the summands are canonical, the splitting is 
certainly not,
as the  choice of the line bundle yields  at least three different splittings  by means of 
three  mechanisms in a  triangulated category introduced by Deligne. It is known that these
three choices yield splittings of mixed Hodge structures.
In this paper, we use the relative hard Lefschetz theorem and elementary linear algebra to  construct
five distinct splittings,
two of which seem to be new, and to prove that they are splittings of mixed Hodge structures.
\end{abstract}

\tableofcontents

\section{Introduction and main theorem}\la{rrivvo}

Let $f: X \to Y$ be a projective map of complex quasi projective varieties,
let $H:=\oplus_{d \geq 0} I\!H^d(X,\rat)$
be the total intersection cohomology rational vector space of $X$
and let $\eta \in H^2(X,\rat)$ be the first Chern class of
an $f$-ample line bundle on $X$. We refer to the survey
\ci{bams} for the background concerning perverse sheaves and the decomposition
and relative hard Lefscehtz theorems
etc. that we use. 

The map $f$ endows $H$ with the perverse Leray filtration $P$. 
The graded objects  $H_p:= P_{p}H/P_{p-1}H$
are non-trivial only in a certain interval $[-r,r]$, with $r=r(f) \in \zed_{\geq 0}$. We thus obtain
the  two objects ${\Bbb H}:= (H,P)$ and ${\Bbb H}_*= \oplus_p (H_p, T[-p])$
($T[-p]=$ the trivial filtration translated to position $p$)
in $\mathscr{V}_\rat \mathscr{F}$, the filtered category  of finite dimensional rational
vector spaces. Obviously,  every filtration
on a vector space by vector subspaces
splits and  we have a
good (= inducing the identity on the graded pieces) isomorphisms $\varphi: {\Bbb H}_* \cong {\Bbb H}$ in 
$\mathscr{V}_{\rat}\mathscr{F}$.

The vector space $H$ underlies a natural mixed Hodge structure (MHS).
The subspaces $P_p H$ are mixed Hodge substructures (MHSS)
so that the graded objects $H_p$ are endowed with a natural MHS.
Let ${\mathscr M}HS$ be the Abelian category of rational  mixed
Hodge structures. It is natural to ask whether there are
good-isomorphisms $\varphi: {\Bbb H}_* \cong {\Bbb H}$ in
$\mathscr{M}HS\mathscr{F}$, the filtered category  of mixed Hodge structures.
In English, do we have splittings $\varphi: \oplus_p H_p \cong H$
such that the component $H_p \to H$ is a map of MHS so that, in particular,
the image is a MHSS?

In this paper, we list five  distinct  such mixed-Hodge theoretic good splittings.
They are built by using the $f$-ample  $\eta$  and they depend on it (see Theorem \ref{tm1}):
\beq\la{eccoli}
\xymatrix{
\omega_{\rm I}(\eta), \;\omega_{\rm II}(\eta),\;
\phi_{\rm I}(\eta), \; \phi_{\rm II}(\eta) , \; \phi_{\rm III}(\eta)
\; : \;\;\;
{\Bbb H}_*  \ar[rr]^{\hskip 2.2cm \cong} &&
{\Bbb H}  \qquad {\rm in} \;\; \mathscr{M}HS\mathscr{F} 
}
\eeq

 The key to our approach is the relative hard Lefschetz theorem (RHL).
 The cup product with  $\eta$ induces
 an arrow $\eta : {\Bbb H}  \to {\Bbb H}[2](1)$ in
 $\mathscr{M}HS\mathscr{F}$. This means that $\eta: H \to H(1)$ (Tate shift $(1)$)
 is such that $\eta: P_p H \to P_{p+2}H(1)$ (translation of filtration $[2]$
 and Tate shift $(1)$).  RHL yields isomorphisms in $\mathscr{M}HS$:
 \beq\la{rrff56}
 \xymatrix{
 \eta^{k} \, : H_{-k} \ar[rr]^{\hskip -0.8cm\cong} && H_{k}(k), \qquad \forall \, k \geq 0.
 }
 \eeq 
 
 Our main technical result, which in fact
 is proved in an elementary way, is as follows:
  let ${\mathscr A}$ be  an Abelian
 category with  shift functors $(n)$ and let   $({\Bbb V},e)$
 be  a pair where  $e: {\Bbb V} \to {\Bbb V} [2](1)$ is an arrow in 
 the filtered  category $\mathscr{AF}$  inducing
 isomorphisms $e^k: V_{-k} \cong V_{k}(k)$, for every $k \geq 0$; then
 there is a natural  isomorphism
 $\omega_{\rm I} (e) : {\Bbb V}_* \cong {\Bbb V}$ in $\mathscr{AF}$.
 
  With this result in hand, we easily verify that
 we can construct the remaining four  splittings within
 $\mathscr{AF}$. We then set $\mathscr{A}= \mathscr{M}HS$
 and deduce (\ref{eccoli}). Let us stress again, that  we use 
 RHL in an essential way and that the point made in this paper is that
 once you have this deep result, the splittings (\ref{eccoli})
 stem from  elementary linear algebra considerations.

 The construction of the three splittings of type $\phi$ is
 borrowed from \ci{delseattle}. However, it seems that
 \ci{delseattle} only yields
  $\phi$-type splittings in $\mathscr{V}_\rat \mathscr{F}$, i.e. not necessarily in its refinement  $\mathscr{M}HS\mathscr{F}$.
 By coupling \ci{delseattle} with the theory of mixed Hodge modules,
 one can indeed prove that the splittings of type $\phi$
 take place in $\mathscr{M}HS\mathscr{F}$. By way of contrast, as pointed out above,
 the constructions of this paper are based on the elementary construction 
 of  $\omega_{\rm I}(e)$ in $\mathscr{AF}$.
 
 The fact that $\phi_{\rm I}(\eta)$ is  mixed-Hodge theoretic
 had been  proved  in \ci{htadt, decpf2} (projective and quasi projective case, respectively)
 by using the  properties of the cup product, of 
 Poincar\'e duality and  geometric descriptions of the perverse Leray filtration
 $P$ on $H$ associated with the map $f$. The proof 
 that  $\phi_{\rm I}(\eta)$ is an isomorphism
 in $\mathscr{M}HS\mathscr{F}$ that we give here
  is different since it does not use the 
 aforementioned special features of the geometric situation.

The simple Examples \ref{snzdiff} and \ref{rrrr} show that
the five splittings (\ref{eccoli}) are, in general, pairwise distinct.

There is a natural condition,
the existence of an $e$-good splitting,
under which the five splittings coincide; see Definition \ref{gspltz}
and  Proposition \ref{eunico}. 

In the paper \ci{dechaumig}, we proved the
following result (auxiliary to the main result of \ci{dechaumig}): the Hitchin fibration
$f: X\to Y$
for the groups $GL_2 (\comp)$, $SL_2 (\comp)$ and $PGL_2(\comp)$
associated with any compact Riemann surface of genus $g \geq 2$ 
and with Higgs bundles of odd degree, presents
 a natural $f$-ample line bundle
$\alpha$ and, in the terminology of the present  paper,
the splitting $\phi_{\rm I} (\alpha)$ is $\alpha$-good
(\ci{dechaumig} shows that (\ref{rfqqlm}) holds for $\phi_{\rm I} (\alpha)$).
In particular, in this case,  the five splittings coincide.

\bigskip
{\bf Acknowledgments.}
I thank Luca Migliorini for useful conversations.
This paper was written during my wonderful stay
at the Department of Mathematics at the University of Michigan,
Ann Arbor, as  
Frederick W. and Lois B.
Gehring Visiting Professor of Mathematics in the Fall of 2011.
I thank the Department of Mathematics of the University of Michigan,
especially Mircea Musta\c{t}\u{a}, for their kind and generous  hospitality
and gratefully acknowledge partial financial support
from the N.S.F., the Frederick W. and Lois B.
Gehring Visiting Professorship fund and the 
David and Lucille Packard foundation.

 \subsection{The main theorem}\la{sstmt}
 A splitting ${\Bbb H}_* \cong {\Bbb H}$ in $\mathscr{V}_{\rat} \mathscr{F}$
 acquires significance only if we can describe $H_*$. 
 This is the content of the decomposition theorem of Beilinson, Bernstein, Deligne and Gabber
 (see the survey \ci{bams}) which implies
 the highly non-trivial fact 
 that, up to a simple renumbering of cohomological degrees,
 we have that $H_p = H(Y, \phix{p}{Rf_* IC_X})$,
 where $\phix{p}{Rf_* IC_X}$ is the $p$-th
 perverse cohomology sheaf of the push-forward $Rf_* IC_X$ of the intersection cohomology complex $IC_X$
 of $X$ ($= \rat_X[\dim_\comp{X}]$, if $X$ is nonsingular).
 
 Let us briefly discuss how the cohomology groups of these perverse sheaves
 split according 
  to the decomposition by supports and to the primitive Lefschetz decomposition coming from RHL.
 Let us start with the one by supports, i.e. the $Z$'s appearing in what follows:
 each perverse sheaf ${\cal H}^p:=\phix{p}{Rf_* IC_X}$ is a semi-simple
 perverse sheaf and  decomposes canonically
 by taking supports: ${\cal H}^p=\oplus_Z {\cal H}^p_Z$, where
 the sum is finite and
 the summands are the intersection cohomology complexes
 with  suitable semisimple local coefficients 
 of suitable closed integral subvarieties  $Z$  of $Y$.
 The
 RHL induces the  primitive Lefschetz
 decompositions (PLD) of
 ${\cal H}^p$:  recall that $p \in [-r,r]$; let $i\geq 0$; set ${\cal P}^{\eta}_{-i,Z} = 
 \ke \, \{ \eta^{i+1}:
 {\cal H}^{-i}_Z \to {\cal H}^{i+2}_Z \}$ and, for $0 \leq j \leq i$,
 set ${\cal P}^{\eta}_{ij,Z}:= \im \, \{\eta^j: {\cal P}^{\eta}_{-i,Z} 
 \to {\cal H}^{-i+2j} \}$; then the PLD reads as
 ${\cal H}^p = \oplus_{-i+2j=p} \oplus_Z {\cal P}^{\eta}_{ij,Z}$ (sum subject
 to $ 0 \leq j \leq i$).
 Set $P_{-i}^{\eta}(-j)_Z:= H(Y, {\cal P}^{\eta}_{ij,Z})$; it is a MHSS of $H_{-i+2j}$.
 Note that $\eta$ induces isomorphisms $\eta^j: (P^\eta_{-i}(0)_Z)(-j) \cong 
 P^{\eta}_{-i}(-j)_Z$ 
 of MHS. 
 By combining the decomposition by supports with the primitive one, we obtain
 the
 splitting $\oplus_{i,j,Z} P^\eta_{-i}(-j)_Z \cong H_*$ which, in view of 
   \ci{htam, htadt, decpf2, absho}, is a splitting
 of MHS.
   We set ${\Bbb P}^\eta_p: = (\oplus_{-i+2j=p, Z} P^\eta_{-i}(-j)_Z,T[-p])$
 (sum subject to $0 \leq j \leq i$; trivial filtration translated to position $p$)
 which is
  an object in $\mathscr{M}HS\mathscr{F}$.
 
  By taking into account these refined 
 splittings, we have the corresponding  refinements
 of (\ref{eccoli}). 
 We note that everything holds just as well, with the same proofs,
  for intersection cohomology with compact supports.

 The main result of this paper is the following
 
 \begin{tm}
 \la{tm1}
 Let $f: X \to Y$, $\eta$  and ${\Bbb H}$ be as above.
 There are five distinguished good splittings  {\rm (\ref{eccoli})} in $\mathscr{M}HS\mathscr{F}$:
 \beq\la{eqtm1}
 \xymatrix{
 \bigoplus_{0\leq j \leq i, \, Z} {\Bbb P}^\eta_{-i}(-j)_Z 
 \ar[rr]^{\hskip 1.2cm \cong} && {\Bbb H}_*.
 }
 \eeq
 Similarely, for intersection cohomology with compact supports.
 \end{tm}
 {\em Proof.} We first work in the abstract setting outlined in
 \S\ref{filcat} 
 of the filtered category $\mathscr{AF}$ of an Abelian category
 $\mathscr{A}$
 endowed with a shift functor.
 In this context, we work with an object ${\Bbb V}$ 
and an arrow $e: {\Bbb V} \to {\Bbb V}[2](1)$ subject to the
HL condition (\ref{bvb}). 

\n
We construct the first splitting
$\omega_{\rm I} (e)$ in Proposition
\ref{splz}, which is the iterated version of the splitting mechanism
established in the elementary Lemma \ref{keylm}. 

\n
The splitting $\omega_{\rm II}(e)$ is obtained
by the ``dual" procedure  as follows.
Let $\mathscr{A}^o$ be the category opposite to ${\mathscr A}$;
it is Abelian and can be endowed with a shift functor
coming from the given one in $\mathscr{A}$.
We  then  have that  $e^o: {\Bbb V}^o \to {\Bbb V}^o [2]^o(1)^o$ in $\mathscr{A}^o\mathscr{F}$
satisifes the corresponding HL condition (\ref{bvb}).
We thus obtain $\omega_{\rm I} (e^o): {\Bbb V}_*^o \cong
{\Bbb V}^o$. We set $\omega_{\rm II}:= ( (\omega_{\rm I}(e^o))^o   )^{-1}.$

\n
The splitting $\phi_{\rm I} (e)$ is constructed in \S\ref{dspfag}.
The proof is parallel to \ci{delseattle}, \S2 with the following
two  changes:
(i) instead of using the existence of a splitting arising from \ci{delseattle}, \S1, which
is proved using some basic features of $t$-categories, features
that $\mathscr{AF}$ does not present, we use the existence 
of either of the splittings $\omega (e)$ established above;
(ii)  we
adapt  the proof of \ci{delseattle}, Lemma 2.1,
which again takes place in the context of a $t$-category,
to the context  of $\mathscr{AF}$.

\n
The splitting $\phi_{\rm II} (e)$ is obtained
in \S\ref{zcdds} by following  the  procedure ``dual"
to the one followed to produce $\phi_{\rm I}(e)$.

\n
Finally, $\phi_{\rm III} (e)$, which necessitates that we work
with a $\rat$-category (= the Hom-sets are rational vector spaces),
is constructed in \S\ref{ztdds} by adapting the corresponding construction
in \ci{delseattle}, \S3. This construction is self-dual.

\n
We now specialize to $\mathscr{A}:= \mathscr{M}HS$, 
with shift functor given by the Tate shift and to $({\Bbb V}, e)
:= ({\Bbb H}, \eta)$ and conclude, due to the fact that
the HL condition (\ref{bvb}) is met
by the RHL (\ref{rrff56}).
\blacksquare

\section{The five splittings}
\la{the5spl}

\subsection{Filtered category associated with an Abelian category}\la{filcat}
Let $\mathscr A$ be an Abelian category whose elements we 
denote $V,W$, etc. 
For ease of exposition only, in this paper
we make heavy use of  the language of sets and elements.

A filtration $F$ on $V$ is a finite increasing filtration, i.e.
an increasing sequence of subobjects $\ldots \subseteq
F_pV \subseteq F_{p+1} V \subseteq \ldots \subseteq $
of $V$ such that $F_pV=0$ for $p \ll 0$ and $F_p V =V$ for $p \gg 0$.
We set ${\rm Gr^F_p V}:= F_pV/F_{i-1}V$. We denote  by $T$ the trivial
filtration  on $V$: $T_{-1}V=0 \subseteq T_0V=V$.
Given $n\in \zed$, we denote by $F[n]$ the $n$-th translate
of $F$: $F[n]_pV:= F_{n+p}V$, so ${\rm Gr}^{F[n]}_pV= {\rm Gr}^F_{n+p} V$;
for example $T[-p]$ is the trivial filtration in position $p$.

Given a pair ${\Bbb V}:= (V,F)$ and a subquotient $U$ of $V$, the filtration $F$ on $V$
induces a filtration on $U$, which we still denote by $F$;
for example  $({\rm Gr}^F_p V,F)=({\rm Gr}^F_p V,T[-p])$. 
In particular, for every $p\leq q \in \zed$, we have the following  pairs
associated with ${\Bbb V}$:
\beq\la{ecvrt}
\xymatrix{
{\Bbb V}_{\leq p} := & (F_p V,F),  &&
{\Bbb V}_{\geq p} := & (V/F_{p-1}V,F),\\
{\Bbb V}_p:= & ({\rm Gr}^F_p V,F), &&
{\Bbb V}_{[p,q]} :=& (F_{q} V/ F_{p-1}V ,F).
}
\eeq
We say that ${\Bbb V}$ has type $[a,b]$, with $a \leq b$, if ${\rm Gr}^F_p V=0$
for every $p \notin [a,b]$.

Given $(V,F)$ and $(W,F)$, an arrow ${\frak l}: V\to W$ in $\mathscr{A}$ is said to be
a filtered arrow if it respects the given filtrations, i.e.
$\frak l$ maps $F_p V $ to $F_pW$, for every $p$. 

The filtered category ${\mathscr AF}$ associated with ${\mathscr A}$
is the category with objects the pairs ${\Bbb V}=(V,F)$
and arrows the filtered arrows. In particular,
an arrow ${\frak l}: {\Bbb V} \to {\Bbb W}$ induces arrows on the objects
listed in (\ref{ecvrt}), e.g. ${\frak l}_i: {\Bbb V}_i \to {\Bbb W}_i$.
An arrow $\frak l$ in $\mathscr{AF}$ is an isomorphism 
if and only if ${\frak l}_i$ is an isomorphism for every $i \in \zed$.

We have   functors
$[n]: \mathscr{AF} \to \mathscr{AF}$, $ (V,F)= {\Bbb V} \mapsto
 {\Bbb V} [n] := (V,F[n])$, etc.

We have the following 
 graded-type  objects, in $\mathscr A$, $\mathscr{AF}$ and $\mathscr{AF}$, respectively:
\beq\la{gdfvv}
V_* \, :=\, \oplus_p V_p,
\qquad 
{\Bbb V}_* \, := \; \bigoplus_p {\Bbb V}_p, \qquad
\widetilde{\Bbb V}_* \, := \; \bigoplus_p \left( V_p, T[-p] \right) .
\eeq

We say that ${\Bbb V}$ splits in
$\mathscr{AF}$ it there is an isomorphism in $\mathscr{AF}$:
\beq\la{vsplits}
\xymatrix{
\varphi \, : \; {\Bbb V}_* \ar[r]^{\hskip 0.4 cm \cong} & {\Bbb V}.
}
\eeq
We say that a splitting $\varphi$ is good if, for every  $p$,  the induced map $\varphi_p:
V_p \cong V_p$ is the identity.

\begin{rmk}\la{precomp}
{\rm
If  ${\Bbb V}$  splits, then there is a good splitting:
let $\varphi_p: {\Bbb V}_p \cong {\Bbb V}_p$
be the induced  isomorphisms and
replace $\varphi$ with $\varphi \circ (\sum_p \varphi_p^{-1})$.
}
\end{rmk}

The category $\mathscr{AF}$ is
pre-Abelian (additive with kernels and cokernels), hence pseudo-Abelian
(every idempotent has a kernel).  
In particular, given an
idempotent $\pi: {\Bbb V} \to {\Bbb V}$, $\pi^2=\pi$,  
we have canonical splittings in $\mathscr{AF}$: 
\beq\la{rfvoo}
{\Bbb V} \; =  \; \ke\,  ({\rm id} -\pi) \oplus \ke\, \pi 
\; =\;  \im\, \pi   \oplus  \ke\, \pi.
\eeq
The arrow 
 $\iota :(V,T) \to (V,T)[1)$  induced by the identity
 is such that 
the induced  arrow 
$\coim \, \iota \to \im \, \iota$ is not an isomorphism
so that $\mathscr{AF}$ is not Abelian.

\begin{ex}\la{exinmd}
{\rm
The example we have in mind is the one where ${\mathscr A}$
is the Abelian category ${\mathscr M}HS$ of integral (or rational) mixed Hodge structures
(MHS)
where the arrows are the maps that respect the weight and Hodge filtrations.
The Tate shift functor, denoted by  $(1)$ is such that
if $\zed$ is the pure Hodge structure with weight zero and type
$(0,0)$,
then $\zed (1)$ is the pure Hodge structure  with  weight $-2$
and type $(-1,-1)$. Note that, for example, the cup product
with the first Chern class  $L$ of a line bundle on a complex algebraic variety
$X$ induces, for every $k\geq 0$, a map $L: I\!H^k(X,\zed) \to I\!H^{k+2}(X,\zed)(1)$.
An element ${\Bbb M}=(M,F)$ of ${\mathscr M}HS \mathscr{F}$ is 
a MHS $M$ (with its weight and Hodge filtrations) equipped with
an {\em additional} filtration $F$ for which $F_p M$ is a mixed Hodge
substructure (MHSS) of $M$ for every $p$.
}
\end{ex}

Let $(1): {\mathscr A} \to {\mathscr A}$, $V \mapsto V(1)$, ${\frak l} \mapsto
{\frak l}(1)$, be an  additive and exact autoequivalence
and, for $m \in \zed$,  denote by $(m)$ its $m$-th iterate, called the $m$-shift
functor. 

By exactness, the  shift functors lift to functors (also called shift functors): 
\beq\la{nnszft}
(m) \,:\, {\mathscr AF} \lorw  {\mathscr
AF}, \quad {\Bbb V}=(V,F)\longmapsto (V(m),F)=: {\Bbb V}(m),
\quad {\frak l} \longmapsto {\frak l} (m),
\eeq
and we have
\beq\la{rvmqqq}
({\Bbb V} (m))_p = {\Bbb V}_p (m), \qquad
{\rm Gr}^F_p ( {\frak l} (m))= ( {\rm Gr}^F_p \,{\frak l})(m).
\eeq
The shift functors commute with the tanslation functors, so that we can write
${\Bbb V}[n](m)$ unambiguously. We have $({\Bbb V}[n](m))_p = {\Bbb V}_{n+p}(m)$, 
etc.

For every $p \in \zed$,
an arrow ${\frak l}: {\Bbb V} \to {\Bbb W}[n](m)$ induces  arrows
in ${\mathscr A}$:
\beq\la{dsds}
\xymatrix{
{\frak l}_p \, : \; {V_p} \ar[r] & W_{n+p}(m),
}
\eeq
where it is understood that: 
\beq\la{rdz}
W_{n+p}(m) \,= \,{\rm Gr}^F_{n+p} (m) \, = \, {\rm Gr}^{F[n]}_p W(m).
\eeq

If ${\Bbb V}$ has type $[a,b]$, then  we have canonical arrows:
\beq\la{rvbt}
\xymatrix{
{\Bbb V}_a \ar[r]^{i_a} & {\Bbb V} \ar[r]^{p_b} & {\Bbb V}_b,
}\eeq
first  inclusion and last quotient,  with compositum $\delta_{ab}{\rm Id}$.

An arrow ${\frak l}[n](m)$ obtained from an arrow ${\frak l}: {\Bbb V} \to {\Bbb W}$
by shift/translation 
is  simply denoted by ${\frak l}:{\Bbb V}[n](m) \to {\Bbb W} [n](m)$ (e.g. the arrow ${\frak l}_a^{-1}$ in (\ref{bbhl})
is really ${\frak l}_a^{-1} (-m)$).

Let ${\Bbb V}$ and ${\Bbb W}$ be in $\mathscr{AF}$, let $m,n \in \zed$,
let ${\frak l}: {\Bbb V}_* \to {\Bbb W}_*[n](m)$ be an arrow in $\mathscr{AF}$
and let ${\frak l}_{pq}: {\Bbb V}_p \to {\Bbb W}_q[n](m)$ be the  $(p,q)$-th component
of ${\frak l}$. We define the degree $d \in \zed$ homogeneous part ${\frak l}^{ \{d \} }$
of ${\frak l}$ by setting:
\beq\la{rfvoqa}
{\frak l}^{ \{d\} } := \sum_{q-p=d} {\frak l}_{pq}, \qquad \qquad  \left({\frak l} = 
\sum_d {\frak l}^{  \{d\} }\right).
\eeq
Since ${\Bbb V}_p = (V_p, T[-p])$ and  ${\Bbb W}_q [n](m) = (W_q(m), T[-q+n])$, we must have
\beq\la{ecvbtqaa}
{\frak l}^{ \{d\}  } =0 , \quad \forall d \geq n+1.
\eeq

Let ${\mathscr A}^o$ denote the Abelian category opposite to ${\mathscr A}$.

\begin{rmk}\la{itsdualstpd}
{\rm
Let $\mathscr{A}=\mathscr{V}_\rat$ be the category of finite dimensional
rational vector spaces and linear maps. The natural contravariant
functor $\mathscr{V}_\rat 
\to \mathscr{V}^o_\rat$ can be identified
with taking dual vector spaces and transposition of linear maps.
Similarely, if we take $\mathscr{A}=\mathscr{M}HS$. 
This  observation may make what follows more down-to-earth
and the computations of explicit examples
easier.
}
\end{rmk}

We have the exact anti-equivalence $(-)^o: {\mathscr A} \to {\mathscr A}^o$,
$(V \stackrel{f}\to W) \longmapsto (V^o  \stackrel{f^o}\leftarrow W^o)$
whose second iterate is the identity functor. 
We endow ${\mathscr A}^o$ with the additive and 
exact shift functors $(m)^o$: $V^o \mapsto V^o (m^o): = (V(-m))^o$.

A filtered object ${\Bbb V} = (V,F)$ in $\mathscr{AF}$ gives rise to a filtered object
${\Bbb V}^o= (V^o, F^o)$ in $\mathscr{A}^o \mathscr{F}$ by setting:
\beq\la{kawer}
F^o_i V^o = (V^o)_{\leq i}: = (V_{\geq -i})^o.
\eeq
Contemplation of the following diagram may be useful:
\beq\la{cotdgr}
\xymatrix{
V_{\leq i-1} \ar[r]^{\rm mono} \ar[dr] & V_{\leq i} \ar[d]  &
(V_{\leq i-1})^o   & (V_{\leq i})^o  \ar[l]_{\rm  \hskip 0.4cm epi} 
 \\
& V  \ar[d] &  
& (V)^o \ar[u] \ar[ul] &
\\
 & V_{\geq i+1} &  V_{\geq i} \ar[l]_{\rm \hskip 0.3cm epi} 
 & (V_{\geq i+1})^o \ar[u]   \ar[r]^{\rm mono} &  (V_{\geq i})^o.
}
\eeq
Clearly,   $(V^o)_i=(V_{-i})^o$ and we set $F^o [n] := (F[-n])^o$.

We obtain  an anti-equivalence $(-)^o:\mathscr{AF} \to \mathscr{A}^o \mathscr{F}$
whose second iterate is the identity functor.
 The anti-equivalence $(-)^o$ is anti-compatible with
translations, shifts and taking graded pieces, etc., for example:
 \beq\la{kkssddee}
 {\Bbb V}^o [n]^o (m)^o \, = \, ({\Bbb V} [-n](-m))^o.
 \eeq
 An arrow ${\frak l} : {\Bbb V} \to {\Bbb W} [n](m)$ in $\mathscr{AF}$
 yields the  arrow ${\frak l}^o: {\Bbb W}^o \to {\Bbb V}^o [n]^o(m)^o$
 in $\mathscr{A}^o\mathscr{F}$.
 This arrow is really 
${\frak l}^o [n]^o (m)^o$, but we omit those decorations for arrows.

\bigskip
We record the following fact for use in the next section.

\begin{lm}\la{quechesve}
Let $\mathscr B$ be an additive category and let $\rho: B \to B'$ be an arrow in $\mathscr B$.
Assume that the kernel  $\iota_{\rho}: \ke\, \rho \to B$ of $\rho$  exists and that there is a
splitting $r: B' \to B$ of $\rho$, i.e. $\rho \circ r= {\rm id}_{B'}$.
Then the natural arrow 
\beq\la{zzx}
\xymatrix{
B' \oplus \ke \, \rho \ar[rr]^{\hskip 0.5cm r + \iota_{\rho}} && B
}
\eeq
is an isomorphism in ${\mathscr B}$.
\end{lm}
{\em Proof.}
Note that $\rho  \circ (1- r \,\rho)=0$, so that there is an unique arrow $s: B \to \ke \,\rho$
such that $1- r \circ  \rho = \iota_{\rho} \circ s$. It is easy to verify that
the arrow:
\beq\la{mmk}
\xymatrix{
B \ar[rr]^{\hskip -0.5cm (\rho, s)} && B' \oplus \ke\, \rho
}
\eeq
yields the desired inverse to $r + \iota_\rho$.
\blacksquare
\begin{rmk}\la{wqa}
{\rm
Assume, in addition,  that  $\mathscr B$ is pseudo-Abelian, e.g. $\mathscr{B=AF}$, and consider the idempotent arrow
 $\pi:= r\circ \rho$.
Then we have a canonical isomorphism  $B=  \im \, \pi  \oplus \ke \, \pi$.
This isomorphism can be canonically identified with the one in (\ref{mmk}), for $\ke \, \pi = \ke \,
\rho$, and $r$ identifies $B'$ with $\im \, \pi$.
}
\end{rmk}

\subsection{A splitting mechanism in $\mathscr{AF}$}\la{splmech}
Let ${\Bbb V}=(V,F)$  in $\mathscr{AF}$ be of type $[a,b]$, 
let $m \in \zed$ and let:  
\beq\la{bbmap}
\xymatrix{
{\frak l}: {\Bbb V}  \ar[r] &
 {\Bbb V} [b-a] (m)
 }
\eeq
be an arrow such that the resulting arrow (\ref{dsds}) is an isomorphism
in $\mathscr{A}$:
\beq\la{bbhl}
\xymatrix{
{\frak l}_{a} : V_{a} \ar[r]^{\hskip -0.1cm \cong} & V_b (m).
}
\eeq

There is the commutative diagram in $\mathscr{AF}$ (see (\ref{rvbt})): 
\beq\la{mmnn}
\xymatrix{
{\Bbb V}_a \oplus {\Bbb V}_b  \ar[rr]^{\hskip 0.3cm r}_{\hskip 0.6cm i_a + {\frak l}i_a {\frak l}_a^{-1}} 
\ar[rrrrrrddd]_{1_{{\Bbb V}_a \oplus {\Bbb V}_b}}& & {\Bbb V}
 \ar[rrrr]^{\rho'}_{({\frak l}_a^{-1} p_b{\frak l}, p_b)} 
\ar[rrrrddd]^{\rho}&&&&
{\Bbb V}_a \oplus {\Bbb V}_b \ar[ddd]_u^{
\left(
\begin{array}{cc}
1_{{\Bbb V}_a} & - {\frak l}_a^{-1} p_b {\frak l}^2 i_a {\frak l}_a^{-1} \\
0 & 1_{{\Bbb V}_b}
\end{array}
\right)
}
\\
&&&&&&\\
&&&&&&\\
&&&&&& {\Bbb V}_a \oplus {\Bbb V}_b
}
\eeq
so that:
\beq\la{qq}
\xymatrix{
\rho   \,: \; {\Bbb V} \ar[rrrrr]^{\hskip -0.4cm \left( {\frak l}_a^{-1} p_b  
{\frak l} \, -\, {\frak l}_a^{-1}  p_b {\frak l}^2 i_a
{\frak l}_a^{-1} p_b \; , \; p_b \right)}&&&&& {\Bbb V}_a \oplus  {\Bbb V}_b.
}
\eeq
The kernel ${\Bbb K}{\rm er}\, {\rho}$ of $\rho$ in $\mathscr{AF}$ is the kernel
$\ke\, \rho$  of the underlying
map in $\mathscr{A}$ with the  filtration induced  by $(V,F)$. The natural 
inclusion induces a map in ${\mathscr AF}$:
\beq\la{mapiota}
\xymatrix{\iota_{\rho}: {\Bbb K}{\rm er}\, {\rho} \ar[rr]&&  {\Bbb V}.
}
\eeq

\begin{rmk}\la{piufaz}
{\rm
Since the arrows $u$ and ${\frak l}_a$ in (\ref{mmnn}) are  isomorphisms, we have that: 
\beq\la{vneor}
\ke \, \rho = \ke \, \rho' = \ke\, (p_b \circ  {\frak l}) \cap \ke \,p_b,
\eeq
and similarely, if we take into account the induced filtrations.
}
\end{rmk}

\begin{lm}\la{keylm} 
The following  arrow is an isomorphism in $\mathscr{AF}$:
\beq\la{fiso}
\xymatrix{ w:
{\Bbb V}_a \oplus {\Bbb K}{\rm er}\, {\rho} \oplus {\Bbb V}_b
 \ar[rrr]^{\hskip 1.4cm i_a +\iota_{\rho} + {\frak l} i_a {\frak l}_a^{-1} }_{\hskip
1.2 cm \cong}&&& {\Bbb V}.
 }
\eeq
\end{lm}
{\em Proof.} Apply Lemma \ref{quechesve}.
\blacksquare

\begin{rmk}\la{stali}
{\rm
 The map $w$ (\ref{fiso}) is uniquely determined by,   and depends on,  ${\frak l}$. However, the component $i_a$, being the inclusion of the first subspace of the
 filtration, is independent of ${\frak l}$.
 }
\end{rmk}

The object ${\Bbb K}{\rm er}\, {\rho}$ is of type $[a+1, b-1]$,
the inclusion $\iota_{\rho}$ induces natural isomorphisms:

\beq\la{rvf}
\xymatrix{ {\iota_{\rho}}_p:  ({\Bbb K}{\rm er}\, {\rho})_p \ar[rr]^{\hskip -1.3cm \cong} && {\Bbb V}_p\,, \qquad \forall \, p \in [a+1, b-1]
}
\eeq
and, by taking subquotients, a natural isomorphism:
\beq\la{kas}
\xymatrix{ {\iota_{\rho}}_{[a+1,b-1]}:
{\Bbb K}{\rm er}\, {\rho} \ar[rr]^{\hskip 0.4 cm \cong} && {\Bbb V}_{[a+1, b-1]}.
}
\eeq
By combining (\ref{fiso}) with (\ref{kas}), we obtain an isomorphism:
\beq\la{vbcn}
\xymatrix{
w_{[a,b]} \, : \;
 {\Bbb V}_a \oplus {\Bbb V}_{[a+1,b-1]} \oplus {\Bbb V}_b 
 \ar[rr]^{\hskip 2.2cm \cong}
 &&
 {\Bbb V},
}
\eeq 
as well as its component:
\beq\la{rrffv00}
\xymatrix{
w_{[a+1, b-1]} \, : \; {\Bbb V}_{[a+1, b-1]}
\ar[r] & {\Bbb V}.
} 
\eeq
Both isomorphisms induce the identity on the $p$-th graded pieces, for every $p$
in (\ref{vbcn}), for $p \in [a+1, b-1]$ in (\ref{rrffv00}).

One may picture the content of Lemma \ref{keylm} as an  un-wrapping of
 the outmost layer ${\Bbb V}_a \oplus {\Bbb V}_b$ of ${\Bbb V}$ via $\frak l$.
 Note that in general, there is no natural non-trivial arrow from
 a subquotient of an object to the object itself. The arrow (\ref{rrffv00}) is made possible
 by the HL condition.

\subsection{The  splittings $\omega_{\rm I} (e)$ and $\omega_{\rm II}(e)$}\la{hlcond}
Let ${\Bbb V}= (V,F)$ in  $\mathscr{AF}$ be  of type $[-r.r]$ for some $r\geq 0$. 
Up to a translation, this condition
can always be met and  leads to   simplified  notation in what follows.

 Let
\beq\la{lev2}
\xymatrix{
e:  {\Bbb V} \ar[rr] && {\Bbb V} [2](1)
}
\eeq
be an arrow in $\mathscr{AF}$. In particular, for every $k\geq 0$,  we have the iterations
$e^k$ and their graded counterparts: (we drop the shift
when denoting a shifted map 
and, in what follows, we drop subscripts for the maps induced on graded objects):
\beq\la{bvb}
\xymatrix{
e^k \,  : \,  {\Bbb V} \ar[r] & {\Bbb V} [2k](k), &
e^k=  V_j \ar[r] & V_{j+2k}(k).
}
\eeq
\begin{ass}\la{hlass}{\rm ({\bf Condition HL}) We assume that $( {\Bbb V}, e)$ satisfies the hard
Lefschetz-type condition (HL), i.e. that
the arrows:
\beq\la{wsx}\xymatrix{
e^k : V_{-k} \ar[rr]^{\hskip 0.2cm\cong} && V_k(k), & \forall k\geq 0,
}
\eeq
are isomorphisms in $\mathscr{A}$.
}
\end{ass}

The following proposition ensures that if the HL condition is met by a given $({\Bbb V},e)$,
then ${\Bbb V}$ splits, i.e. we have an isomorphism ${\Bbb V} \cong {\Bbb V}$
(\ref{vsplits}). 
By keeping with the analogy of Remark \ref{stali}, we may say that HL
allows to completely un-wrap ${\Bbb V}$. 
 Recall the notion
  of  good splitting (on inducing the identity on the graded pieces).

\begin{pr}\la{splz}
Let $( {\Bbb V}, e)$ satisfy the HL condition.
There is a good splitting:
\beq\la{isosplz}
\xymatrix{ \omega_{\rm I} (e) = \omega_{\rm I} \, : \;
{\Bbb V}_*=
\bigoplus_{p \in [-r.r]} {\Bbb V}_p \ar[rr]^{\hskip 2.5cm \cong} && {\Bbb V}.
}
\eeq
\end{pr}
{\em Proof.}
By applying Lemma \ref{keylm} and  to ${\frak l}:= e^r$,  so that $[a,b]=[-r,r]$
and $m=r$,
we obtain
\beq\la{vbvb}
\xymatrix{
w: {\Bbb V}_{-r} \oplus {\Bbb K}{\rm er} \, \rho \oplus {\Bbb V}_r \ar[rr]^{ \hskip 1.8cm \cong} &&
{\Bbb V}.
}
\eeq
The arrow $\frak l$ yields the arrow $\tilde{\frak l}:= \omega_{\rm I}^{-1}\circ 
{\frak l}\circ \omega_{\rm I}$
on the lhs of (\ref{vbvb}).
Keeping in mind that (\ref{fiso}) means that the filtration $F$ on $V$
splits, we obtain 
the ${\Bbb K}{\rm er}\, \rho$-component ${\frak l}'$ of ${\tilde{\frak  l}}$:
\beq\la{kaqwe}
\xymatrix{
{\frak l}' :  {\Bbb K}{\rm er}\, \rho \ar[rr] && {\Bbb K}{\rm er}\, \rho [2](1).
}
\eeq
In view of (\ref{rvf}), we have that ${\frak l}'$
satisfies the HL condition. 

\n
By using (\ref{kas}) and (\ref{vbcn}), we replace  $ {\Bbb K}{\rm er}\, \rho$ with ${\Bbb V}_{[-r+1, r-1]}$ and we obtain  the desired splitting $\omega_{\rm I}(e)$
by descending induction on $r$.

\n
By construction (i.e. $\tilde{\frak l}= \omega_{\rm I}^{-1} \circ  {\frak l}
\circ \omega_{\rm I}$,
(\ref{rvf}) and (\ref{wsx})), the isomorphism $\omega_{\rm I}(e)$ 
 induces the identity on the graded pieces and is thus good.
\blacksquare

\bigskip
The splitting $\omega_{\rm II}(e)$ is obtained by the dual construction.
This is explained in the proof of Theorem \ref{tm1}.

\begin{rmk}\la{nozdr}
{\rm
In general, $\omega_{\rm I} (e) \neq \omega_{\rm II} (e)$; see
Examples \ref{snzdiff} and \ref{rrrr}.
In particular, neither of the two constructions $\omega_{\rm I}(e)$
and $\omega_{\rm II}(e)$ is self-dual.
}
\end{rmk}

\bigskip
The purpose of the next three sections  is to show that if
${\Bbb V}$ splits in $\mathscr{AF}$,  then one can use the HL property
(\ref{wsx}) to construct
three additional natural splittings taking place in $\mathscr{AF}$.
Note that these constructions are based solely on the existence of a  splitting.

\subsection{The first Deligne splitting $\phi_{\rm I}(e)$}\la{dspfag}
Let $({\Bbb V},e)$ be as in (\ref{lev2}) and assume that it  satisfies the HL condition
  (\ref{wsx}). In particular, in view of  (\ref{isosplz}),    ${\Bbb V}$ splits
  in $\mathscr{AF}$.

 Let $i \geq 0$ and 
 define the primitive objects in $\mathscr A$:
 \beq\la{defprim}
 P_{-i} := \ke \, \left\{ e^{i+1}: V_{-i} \lorw V_{i+2}(i+1) \right\}.
 \eeq
The subquotient  $P_{-i}$ of $V$ inherits the filtration induced by
 $F$ on $V$, i.e. the trivial filtration translated in position $-i$, and we denote the resulting object in $\mathscr{AF}$
 by: 
 \beq\la{rgbnyqaa}
 {\Bbb P}_{-i}=(P_{-i}, F)=(P_{-i},T[i]).
 \eeq
 We have the  natural  monomorphisms in ${\mathscr A}$:
 \beq\la{rvqzew}
 \xymatrix{
 P_{-i}(-j)  \ar[r] & V_{-i}(-j) \ar[r]^{\hskip 0.2cm e^j} & V_{-i+2j}
 }
 \eeq
 which, taken together, yield
 the canonical  primitive Lefschetz decompositions (PLD) in
 $\mathscr A$ and $\mathscr{AF}$, respectively:
 \beq\la{pldf}
\xymatrix{
\bigoplus_{0\leq j \leq i}   P_{-i}(-j) \ar[r]^{\hskip 1cm \e}_{\hskip 1cm \cong} &
 {V}_* ,  & \\
 \bigoplus_{0\leq j \leq i}    \left(P_{-i}(-j), T[i-2j]
 \right) =
  \bigoplus_{0\leq j \leq i}   {\Bbb P}_{-i}(-j) [-2j]
  \ar[rr]^{\hskip 4.4cm \e}_{\hskip 4.4cm \cong} &&
 {\Bbb V}_*.
 }
 \eeq

We have the commutative diagram
 of epimorphisms and monomorphisms in $\mathscr{AF}$:
 \beq\la{epiz}
 \xymatrix{&&& {\Bbb P}_{-i} \ar[r]^{\rm mono} \ar[d]_{\rm mono} &{\Bbb V}_{-i} \ar[dl]^{\rm mono}\\
 {\Bbb V}\ar[r]^{\rm \hskip -0.2cm epi} &  \ldots \ar[r]^{\rm 
 \hskip -0.4cm epi} &  {\Bbb V}_{\geq -i-1} \ar[r]^{\rm \hskip 0.1cm epi} & {\Bbb V}_{\geq -i}. &
 }\eeq

 Note that (\ref{rgbnyqaa}) implies that
  if ${\frak l}: {\Bbb P}_{-i} \to {\Bbb W}$ is an arrow in $\mathscr{AF}$,
 then it factors through ${\Bbb W}_{\leq -i}$, i.e. the underlying arrow ${\frak l}:P_{-i} \to W$, factors through $F_{-i}W$.
 
 \begin{lm}
 \la{delsplgr} Let $({\Bbb V},e)$ be as above.
 There is a unique arrow  $f_i: {\Bbb P}_{-i} \to {\Bbb V}$  in $\mathscr{AF}$  with the following properties:
 \ben
 \item
 it lifts the natural arrow ${\Bbb P}_{-i} \to {\Bbb V}_{\geq -i}$ in {\rm (\ref{epiz})};
  \item
 for every $s>i\geq 0$, the  composition of the arrows below is zero:
 \beq\la{iszew}\xymatrix{
  {\Bbb P}_{-i} \ar[rr]^{e^s \circ f_i} && {\Bbb V} (s) \ar[r] & {\Bbb V}_{\geq s} (s).
 }
 \eeq
 \een
 \end{lm}
 {\em Proof.} The proof is essentially identical to the one of
 \ci{delseattle}, Lemme 2.1 (see also \ci{delseattle}, 2.3). We include it,
 with the necessary changes, for the reader's convenience.
 Recall that we use the language of sets. 
 
 \medskip
 \n
 Let $\Phi: \mathscr{AF} \to \mathscr{B}$ be an additive functor into an Abelian category
${\mathscr B}$.
We denote $\Phi(e)$ simply by $e$. 
 Let $i \geq 0$ and $x \in \Phi (V_{\geq -i})$
be such that $0= e^{i+1} (x) \in \Phi( V_{\geq i+2})$.

\medskip
 \n
 {\em CLAIM 1:}
there is a unique lift
$y \in \Phi(V_{\geq -i-1})$ of $x$ such that $0= e^{i+1} (y) \in   \Phi( V_{\geq i+1})$.

\medskip
\n
{\em Proof.} 
For every $a \in \zed$, we have the natural maps: 
\beq\la{edcrtq}
\xymatrix{
{\Bbb V}_a \ar[r] & {\Bbb V}_{\geq a+1} \ar[r] & {\Bbb V}_{\geq a+1}.
}\eeq
Since $\Phi$ is additive and, in view of Lemma \ref{splz}, ${\Bbb V}$ splits in $\mathscr{AF}$, 
we have the  short exact sequences 
in $\mathscr B$ stemming from (\ref{edcrtq}): 
\[
0 \lorw \Phi( {\Bbb V}_a) \lorw \Phi({\Bbb V}_{\geq a})  \lorw \Phi({\Bbb V}_{\geq a+1}) \lorw 0.\]
By naturality, we have the following commutative diagram of short exact sequences
in ${\mathscr B}$: 
\beq\la{beltr}
\xymatrix{
 0 \ar[r] & \Phi ({\Bbb V}_{-i-1}) \ar[r] \ar[d]^{e^{i+1}}_{\cong} & \Phi ({\Bbb V}_{\geq -i-1}) \ar[r] \ar[d]^{e^{i+1}} &
 \Phi ({\Bbb V}_{\geq -i}) \ar[r] \ar[d]^{e^{i+1}} &  0 \\
 0 \ar[r] & \Phi ({\Bbb V}_{i+1}(i+1)) \ar[r]  & \Phi ({\Bbb V}_{\geq i+1}(i+1)) \ar[r] &
 \Phi({\Bbb V}_{\geq i} (i+1)) \ar[r] & 0.
}
\eeq
The claim follows from a simple diagram-chase starting at $x \in \Phi ({\Bbb V}_{\geq -i})$.
 
 \medskip
 \n
{\em  CLAIM 2:} under the same hypotheses as the ones of CLAIM 1, there is a unique lift
 $y \in \Phi({\Bbb V})$  of $x\in \Phi ({\Bbb V}_{\geq i})$ such that $\forall s >i$, we have that
$0 = e^s (y) \in \Phi ({\Bbb V}_{\geq s})$.
 
 \n
 {\em Proof.} The element $y$ found in CLAIM 1 satisfies the hypotheses
 of CLAIM 1 for $i+1$. Since the filtration is finite, we conclude by a repeated use of 
 CLAIM 1.
 
 \medskip
\n
Let us apply CLAIM 2 to the functor  $\Phi(-):= {\rm Hom}_{\mathscr{AF}}({\Bbb P}_{-i}, -): \mathscr{AF}
\to {\rm Ab}$ (Abelian groups): set  $x: {\Bbb P}_{-i} \to
{\Bbb V}_{\geq -i}$  to be  as in (\ref{epiz}).

\n The statement of the lemma follows by applying CLAIM 2 to $x$:
the hypotheses of CLAIM 2 are met in view of the defining property (\ref{defprim}) of $P_{-i}$,
and the resulting element $y$ is the desired $f_i$.
 \blacksquare

 \bigskip
 Let $\varphi: {\Bbb V}_* \cong {\Bbb V}$ be any splitting.
 In view of
 the primitive Lefschetz decomposition (\ref{pldf}), we
 can talk about the components $\varphi_{ij}: {\Bbb P}_{-i}(-j) [-2j] \to {\Bbb V}$.
 Of course, 
  there are many splittings
having components $\varphi_{ii} = f_i$. 

We define the 
{\em first Deligne isomorphism $\phi_{\rm I} (e)$
associated with $( {\Bbb V}, e)$} by taking the compositum
of the following two isomorphism
\beq\la{del1}
\xymatrix{
\phi_{\rm I} (e) \, = \phi_{\rm I} \, : \; 
{\Bbb V}_* \ar[rr]^{\hskip -0.2cm \e^{-1}}_{\hskip -0.4cm (\ref{pldf})} && \bigoplus_{0\leq j \leq i} {\Bbb P}_{-i}(-j) [-2j] \ar[rr]^{\hskip 1.8cm \sum e^j \circ f_i}_{\hskip 1.8cm \cong} &&
{\Bbb V}.
}
\eeq

By using the language of elements, if we  denote by  $p_{ij}$
the typical element in $ P_{-i}(-j)$  and we form
the typical element $v_* \in V_*$:
\beq\la{wsppqq}
v_*\,=\, \sum_p v_p \,=\,  \e \left( \sum_p \sum_{-i+2j=p} p_{ij}
\right)  \;\in \; V_*,
\eeq
 then
\beq\la{evmoq}
 \phi_{\rm I}(e) \, : \; v_* \, \longmapsto  \,   \sum_{0 \leq j \leq i}
 e^{j} \left(f_{i} (p_{ij}) \right).
 \eeq
 In particular, we have that
 \beq\la{ddooqq1}
 \xymatrix{
  \phi_{\rm I} (\e (e^{l} p_{ij})) = e^l \phi_{\rm I} (\e ( p_{ij})), 
  \qquad \forall \, 0 \leq l \leq i-j.
  }
 \eeq
 
 \begin{rmk}\la{tunits}{\rm
Let us omit the shifts, translations and filtrations.
Let $\varphi: {\Bbb V}_* \cong {\Bbb V}$ be a  splitting and $\varphi_i: P_{-i}\to V$
be the resulting components.
By (\ref{defprim}), we have that $e^{i+t} \circ \varphi_i: P_{-i} \to V_{\leq i+2t-1}$,
for every $t>0$. Lemma
 \ref{delsplgr} yields    $f_i: P_{-i} \to V$ with
  $e^{i+t} \circ f_i:
P_{-i} \to V_{\leq i+t-1}$ for every $t>0$, i.e.  an improvement by $t$ units
with respect to an arbitrary splitting, even a good one.
The paper \ci{dechaumig} exploits this special property of $\phi_{\rm I}(e)$
in the context of a  study of the geometry of the
Hitchin fibration.
}
\end{rmk}
 
 \begin{rmk}\la{rvbblmlm}
 {\rm By constuction (see (\ref{del1}) and Remark \ref{precomp}), 
 the isomorphisms $\phi_{\rm I} (e)$  and $\omega_{\rm I}(e)$ are good (\ref{vsplits}).
 In general, the two differ from each other; see Examples \ref{snzdiff} and \ref{rrrr}. However, they
agree on ${\Bbb V}_{-r} \oplus {\Bbb V}_r$: 
in fact, both induce the identity on the graded pieces, so that
they must agree on ${\Bbb V}_{-r}= {\Bbb P}_{-r}$;  by comparing
the expression  (\ref{del1}) for $\phi_{\rm I}(e)$  restricted to ${\Bbb V}_{r}$
with the corresponding one for $\omega_{\rm I} (e)$, i.e. (\ref{vbvb}) and 
(\ref{fiso}), we see that they coincide. In particular, if $V_i =0$
for every $|i|\neq r$, then $\omega_{\rm I}(e) = \phi_{\rm I} (e)$.
}
\end{rmk}

%Any two isomorphism  $\varphi_1, \varphi_2: {\Bbb V}_* \cong {\Bbb V}$
% yield
%$\phi : = \varphi_2^{-1}\circ \varphi_1$ with
%homogeneous parts   $\phi^{ \{ d \} } =0$
%for every $d>0$ (\ref{ecvbtqaa}). If $\varphi_1$ and $\varphi_2$ are both good,
%then $\phi$ is the identity modulo negative degree:
%\beq\la{przzaw}
%\phi = {\rm Id}_{ {\Bbb V}_* } + \sum_{d<0} \phi^{ \{d \}  }, \qquad 
%\phi \equiv {\rm Id}_{ {\Bbb V}_* } \quad  \mbox{(modulo negative degree)}.
%\eeq

Let  $\varphi: {\Bbb V}_* \cong {\Bbb V}$ 
be a splitting. The matrix $\tilde{e}(\varphi)$  of $e$
with respect to $\varphi$ is defined by setting:
\beq\la{rvpqcr}
  \tilde{e} (\varphi)\, = \, \tilde{e} := 
  \left( \varphi^{-1} \circ e \circ \varphi \right) \, : \; 
{\Bbb V}_* \lorw {\Bbb V}_* [2](1),
\qquad
\tilde{e}  = \sum_{pq} \tilde{e}_{pq} = \sum {\tilde{e}}^{ \{d \} }.
\eeq
By virtue of (\ref{ecvbtqaa}), we have that
\beq\la{rrttnnyu}
\tilde{e}^{ \{d\} } =0 , \quad \forall d >2.
\eeq
Let us assume that $\varphi$ is good. Then
\beq\la{ebbuono}
\qquad {\tilde{e}}^{ \{ 2 \}} = \sum_p e, \quad e : {\Bbb V}_p \stackrel{e}\lorw 
{\Bbb V}_{p+2}(1). 
\eeq
Note that while 
$\tilde{e}^{ \{2\} }$
is independent of $\varphi$, we have that
$\tilde{e} (\varphi)^{ \{d\} }$ depends on $\varphi$ for $d \leq 1$.

We have the refinement  $\hat{e}(\varphi)$ of the matrix $\tilde{e}(\varphi)$
that takes into account the primitive Lefschetz decomposition (\ref{pldf}):
\beq\la{vbapqrew}
\hat{e} (\varphi) \, = \, \hat{e} \, : =\;
\left(  \e^{-1} \circ \tilde{e} \circ \e \right) \, : \;
\bigoplus_{0 \leq j \leq i} 
{\Bbb P}_{-i}(-j) [-2j] \lorw 
\bigoplus_{0 \leq j \leq i} 
{\Bbb P}_{-i}(-j) [-2j] [2](1).
\eeq
By taking components, we have arrows
\beq\la{ijarr}
{\hat{e}}(\varphi)_{ij}^{kl} \, : \; 
{\Bbb P}_{-i}(-j) [-2j] 
\lorw {\Bbb P}_{-k}(-l) [-2l] [2](1).
\eeq

Proposition 2.7 in \ci{delseattle} can be easily adapted to the present context
and yield the following characterization of $\phi_{\rm I} (e)$.
For a ``visual", see  \ci{delseattle}, p.119.

\begin{lm}\la{chardel1}
The splitting $\phi_{\rm I}: {\Bbb V}_* \cong {\Bbb V}$
is characterized among the good ones 
by the following two conditions: 
\ben
\item
for $0 \leq j <i$, we have 
${\hat{e}} (  \phi_{\rm I}  )_{ij}^{kl}=0$ except for ${\hat{e}}
(\phi_{\rm I} )_{ij}^{i,j+1} ={\rm Id}$;
\item
for $j=i$, we have  the ${\hat{e}}(  \phi_{\rm I} )_{ii}^{kl} = 0$ 
except, possibly, for $l \leq i$.
\een
\end{lm}

\begin{defi}\la{gspltz}
{\rm 
We say that a splitting $\varphi: {\Bbb V}_* \cong {\Bbb V}$ is {\em $e$-good}
if it induces the identity on the graded pieces  and
$\tilde{e}(\varphi)$ is homogeneous of degree two:
\beq\la{efgtrwq}
\tilde{e} (\varphi ) = \tilde{e} (\varphi )^{ \{2 \} }.
\eeq 
}
\end{defi}

Clearly, $\varphi$ is $e$-good if and only we have that
\beq\la{rfqqlm}
\hat{e}(\varphi)_{ij}^{kl}=0 \qquad \mbox{except for} \;\; 
\hat{e}(\varphi)_{ij}^{i,j+1}= {\rm Id}, \;\forall \,0\leq j\leq i-1.
\eeq
We also have that $\varphi$ is $e$-good if and only
the composita:
\beq\la{lalskdjf}
\xymatrix{
{\Bbb P}_{-i} \ar[r]^\subseteq & {\Bbb V} \ar[rr]^{\hskip -1.4cm e^{i+1}} && {\Bbb V}[2i+2](i+1)
}
\eeq
are zero for every $i\geq 0$. In this case we say that the  $i$-th graded primitive objects
$P_{-i}$ are embedded into $V$ via $\varphi$ as bona-fide 
$i$-th  primitive
classes: i.e. killed  by $e^{i+1}: P_{-i} \to V_{\leq i+1}(i+1)$,  and not just killed by the
subsequent   projection to $V_{i+1} (i+1)$.

\begin{rmk}\la{fgtrrt}{\rm
Lemma \ref{chardel1} implies that
 if $\varphi$ is $e$-good, then
$\varphi=  \phi_{\rm I} (e)$. In particular, if there exists an $e$-good
splitting, then it is unique.
However, $e$-good splittings do not exist in general:
the reader can verify this in Examples \ref{snzdiff} and \ref{rrrr};
 in the latter example, one can even take
$\pn{1} \times \pn{1} \to \pn{1}$.
Proposition \ref{eunico} shows that the existence of
an $e$-good splitting is rare.
}\end{rmk}

\subsection{The second Deligne splitting $\phi_{\rm II} (e)$}\la{zcdds}
Let $({\Bbb V},e)$ be as in  the beginning of
\S\ref{dspfag}.  In particular,   ${\Bbb V}$
  admits  a splitting in $\mathscr{AF}$ as in  Proposition
\ref{splz}: in fact, we have three so far
  $\omega_{\rm I}(e)$, $\omega_{\rm II} (e)$ and $\phi_{\rm I}(e)$.

The first Deligne isomorphism $\phi_{\rm I} (e^o)$ (\ref{del1}) associated with $({\Bbb V}^o, e^o)$
in $\mathscr{A}^o \mathscr{F}$ yields, by application of $(-)^o:
\mathscr{A}^o\mathscr{F} \to \mathscr{AF}$, the isomorphism in $\mathscr{AF}$:
\beq\la{del2pre}
\xymatrix{
(\phi_{\rm I}(e^o))^o\, : \,
{\Bbb V} \ar[rr]^{\hskip 1.0 cm \cong} &&
  {\Bbb V}_*.
  }
  \eeq
  We define the {\em second Deligne isomorphism associated with $({\Bbb V},e)$} to be
\beq\la{del2}
\xymatrix{
\phi_{\rm II}(e)\, = \, \phi_{\rm II}  \, := \;  ((\phi_{\rm I} (e^o))^o)^{-1} \, : \;
{\Bbb V}_*  \ar[rr]^{\hskip 3.0 cm \cong} && {\Bbb V}.
}
\eeq

In this context, the analogue of Lemma \ref{delsplgr}  reads as follows.
Let \beq\la{rfvoiwxc}
\xymatrix{
f_{ij}' \, : \; {\Bbb V} \ar[rr] & & \;{\Bbb P}_{-i}(-j)[-2j]
\;\; (\stackrel{\e}\subseteq {\Bbb V}_i)
}
\eeq
be the  components of (\ref{del2pre})
associated with (\ref{pldf}) ($\e$ as in (\ref{pldf})).
\begin{lm}\la{fvmoooo}
For every $i \geq 0$, 
the arrow $f_{ii}'$ is the unique arrow ${\Bbb V} \to {\Bbb V}_i$ such that:
\ben
\item by taking the $i$-th graded pieces, $f_{ii}'$ induces the natural projection 
${\Bbb V}_i \to {\Bbb P}_{-i}(-i)[-i]$;
\item
for every $s >i$, the composition below is zero  is (see {\rm \ci{delseattle}, \S3.1}):
\beq\la{ewsce}
\xymatrix{
{\Bbb V}_{\leq -s} \ar[r]^\subseteq & {\Bbb V} \ar[r]^{\hskip -0.4cm \e^s} & {\Bbb V}[2s](s) \ar[rr]^{\hskip -0.8cm f_{ii}'} && {\Bbb P}_{-i}(-i)[-2i][2s](s).
}
\eeq

\een
\end{lm}

By using Lemma \ref{fvmoooo}  and the explicit formula   (\ref{evmoq}) for $\phi_{\rm I}(e)$,
it is easy  to deduce the following explicit expression
for the arrows $f_{ii}'$:
\beq\la{rmer}
 f_{ii}' (\phi_{\rm I} (v_*)) = p_{ii} \qquad\mbox{ ($v_* \, = \, \sum_{0\leq j'\leq i'} p_{i'j'}$). }
\eeq

In general, $\phi_{\rm I} (e) \neq \phi_{\rm II}(e)$ and this discrepancy
is due to the fact that  $f_{ij}' (\phi_{\rm I} (v_*) ) \neq  p_{ij}$.
We now discuss how this discrepancy is  measured exactly in terms  the matrix  $\hat{e} (\phi_{\rm I})$
(\ref{ijarr}) of $e$.

By combining (\ref{del1}) with (\ref{del2pre}), we see that
there is the commutative diagram (the bottom identification is due to (\ref{pldf})):
\beq\la{ffggooiq}
\xymatrix{
{\Bbb V} \ar[rrr]^{e^{i-j}} \ar[d]^{f'_{ij}}  &&&{\Bbb V} [2(i-j)](i-j) \ar[d]^{f'_{ii}}\\
{\Bbb P}_{-i}(-j) [-2j] \ar[rrr]^{\hskip -0.4cm e^{i-j}={\rm Id}} &&& 
{\Bbb P}_{-i}(-i) [-2i] [2i-2j] (i-j).
}
\eeq
We fix $0\leq j\leq i$. For every $0\leq s \leq t$, we use 
(\ref{evmoq}) and (\ref{ddooqq1}) together with  (\ref{rmer}) and  (\ref{ffggooiq})
with the goal of  determining the value of
\beq\la{dkfjgh}
f'_{ij}  \left(  e^s f_{t} (p_{ts}) \right).
\eeq
Recalling that (\ref{del2pre}) induces the identity on the graded pieces,
we deduce that:
\ben
\item
if $t=i$ and $s<j$, then $f'_{ij}  \left( e^s f_{t} (p_{ts}) \right)=0$;

we can see this, as well as the assertions that follow, on the following diagram (we do not write $\e$):
\beq\la{3030}
e^s f_t (p_{is})=\phi_{\rm I} (p_{is}) 
\longmapsto e^{i-j} \phi_{\rm I} ( p_{is}) = \phi_{\rm I} (e^{i-j} p_{is})
\stackrel{ f_{ii}'}\longmapsto 0;
\eeq
\item
if $t=i$ and $s=j$ then $f'_{ij}  \left( e^j f_{i} (p_{ij}) \right)= p_{ij}$;
\item
if $t\neq i$ and $s+i -j \leq t$, then  $f'_{ij}  \left(
e^s f_{t} (p_{ts})
\right)=0$;
\item if $t\neq i$ and $\s:=s+i -j -t \geq 1 $, then  
\beq\la{rfvvfrre}
f'_{ij}  \left(e^s f_t (p_{ts})\right)= f'_{ii} (   e^\s e^t f_t (p_{ts})     )
\eeq
which, recalling the definition (\ref{ijarr}) of $\hat{e}_{\rm I} = \hat{e}
(\phi_{\rm I} (e))$, has the following form:
\beq\la{rrii11}
 \qquad
(\hat{e}_{\rm I}^\s)_{tt}^{ii}  \,(q_{tt})
\qquad
\left(\mbox{where $q_{tt} \, := \, f'_{tt} (e^t f_t (p_{ts})$}\right).
\eeq
\een

\begin{pr}
\la{eecc11}
The first Deligne isomorphism $\phi_{\rm I}(e)$ is 
$e$-good {\rm (\ref{rfqqlm})} if and only if
$\phi_{\rm I} (e) = \phi_{\rm II} (e)$.
\end{pr}
{\em Proof.} In view of (\ref{del2pre}) and \ref{del2}),
we have that  the two  Deligne isomorphisms coincide if and only if,
using the notation in 
(\ref{rmer}), we have that
$f'_{ij}( v_*) = p_{ij}$, for every $0 \leq j \leq i$. According to the four points above,
the only obstruction to having this latter condition
stems from (\ref{rrii11}) not being zero for some pair $(i,t)$ with $i\neq t$.
By reasons of degree, i.e. by (\ref{ecvbtqaa}), since $i\neq t$,
we have that
 ${({\hat e}_{\rm I})}^{ii}_{tt}$
 is in degrees $\leq 1$ .  If $\phi_{\rm I}$ is $e$-good, then
${\hat{e}}_{\rm I}$ is of pure homogeneous degree $2$, so that 
 $e^\s e^t f_t (p_{ts}))=0$ and we infer the desired equality. 

\n
Conversely, let us assume that the two Deligne isomorphisms coincide. By contradiction
let us assume that $\phi_{\rm I}$ is not  $e$-good. 
According to (\ref{rfqqlm})   there are integers $0 \leq t$ and $0 \leq l \leq k$
and  a non zero arrow
 $l \leq t$.  Among these non zero arrows ${ (\hat{e}_{I})}_{tt}^{kl}$,
 chose one, ${ (\hat{e}_{I}) }_{t_ot_o}^{k_ol_o}$, for which  the difference $k-l$ attains the minimum value. 
 In the language of elements, what above ensures that  there is $0\neq p_{t_o} \in P_{-t_o}$
 such that:
 \beq\la{eedd44}
 e^{t_o+1} f_{t_o} (p_{t_o}) = e^{l_o} f_{k_o} \left(  \hat{e}_{t_o t_o}^{k_ol_o} 
 \left(  e^{t_o} p_{t_o} \right)   \right)
 + {\sum}^*e^{l} f_{k} \left( \hat{e}_{t_ot_o}^{kl} \left(e^{t_o} p_{t_o} \right)    \right),
 \eeq
where  the first term on the r.h.s. is non-zero and
 ${\sum}^*$ is the sum over the non-zero terms with $(k,l) \neq (k_o,l_o)$.
Since for these latter terms, $k-l \geq k_o-l_o$, we deduce that
\beq\la{ee11}
f'_{k_ok_o} \left(  e^{k_o -l_o}  {\sum}^* \right)  = 0 \in P_{-k_o} (-k_o).
\eeq
On the other hand,  since obviously  $e^{k_o-l_o} e^{l_o}= e^{k_o}$, 
we have that:
\beq\la{ee44}
q_{k_o} \, : = \;
f'_{k_ok_o} \left(  e^{k_o -l_o}    e^{l_o} f_{k_o} \left(  \hat{e}_{t_o t_o}^{k_ol_o} 
 \left(  e^{t_o} p_{t_o} \right)   \right) \right)  \neq 0 \in P_{-k_o}(-k_o).
\eeq
In view of (\ref{ffggooiq}), we have that
\beq\la{113355}
f'_{k_o, l_o-1}  \left( e^{t_o} f_{t_o} (p_{t_o}) \right) = q_{t_o} \neq 0.
\eeq 
Since we are assuming that the two Deligne isomorphisms coincide,
by virtue of the first paragraph of this proof, we must have
$k_o=t_o$ and $t_o = l_o-1$. This contradicts $l_o \leq t_o$.
\blacksquare

\begin{rmk}\la{rfedserwq}
{\rm In general, there is no $e$-good splitting; see 
Examples \ref{snzdiff} and \ref{rrrr}. In particular, 
neither of the two  constructions  
$\phi_{\rm I}(e)$ and $\phi_{\rm II} (e)$ is  self-dual.
}\end{rmk}

\subsection{The third Deligne splitting $\phi_{\rm III}(e)$}\la{ztdds}
Let $({\Bbb V} ,e)$ be as in \S\ref{dspfag}.  In particular,   ${\Bbb V}$
  admits  a splitting (\ref{isosplz}) in $\mathscr{AF}$ as in  Proposition
\ref{splz}; in fact, we have four, so far.
We also  assume that the Abelian category $\mathscr A$ is $\rat$-linear,
i.e. that  ${\rm Hom}$-groups  are rational vector spaces. The reason for this
is that, in what follows,  one needs to exploit the  ${\rm sl}_2 (\rat)$-action  arising from the given arrow  $e$.

The goal of this section  is to construct the third Deligne isomorphism associated with $({\Bbb V},e)$. Whereas we omit
the  detailed presentation of the algebra underlying this construction (see 
\ci{delseattle}, Lemme 3.3 and Proposition 3.5), we  review some of the  key points, state
its characterization and, along the way, indicate the necessary changes.

Let ${\Bbb V}$ and ${\Bbb W}$ be in $\mathscr{AF}$ and set:
\beq\la{troe}
L^{[n]}_{(i,j)} ( {\Bbb V}, {\Bbb W}) : = 
{\rm Hom}_{\mathscr{AF}} ({\Bbb V}(i), {\Bbb W}(j)[n]).
\eeq
Up to  the canonical isomorphism induced by the shift functors,
the above depends only on the difference $m:=(j-i)$ and  we denote the resulting 
bi-functor by
$L^{[n]}_{(m)}$.

Recalling the definition of the graded-type objects (\ref{gdfvv})
and of degree of maps (\ref{rfvoqa})
(an arrow $(V_p, T) \to (V_q,T) [n](m)$ in $\mathscr{AF}$ has degree $d:= q-p$), we have the natural decomposition by homogeneous degrees:
\beq\la{dkdk}
L^{[n]}_{(m)} ( \tilde{\Bbb V}_*,  \tilde{\Bbb V}_* ) = \bigoplus_{d=-2r}^{2r}
L^{[n], \{d\}}_{(m)} ( \tilde{\Bbb V}_*, \tilde{\Bbb V}_*).
\eeq
The arrow $h: = \sum_p p\, {\rm Id}_{(V_p, T)}$ induces
the   arrow:
\beq\la{efbgt}
h: L^{[n]}_{(m)} ( \tilde{\Bbb V}_*,\tilde{\Bbb V}_*) \lorw L^{[n]}_{(m)} ( \tilde{\Bbb V}_*, \tilde{\Bbb V}_*),
\quad
u \longmapsto h\circ  u;
\eeq
this arrow is of homogeneous degree zero, i.e. $\{d\} \mapsto \{d\}$, with respect to (\ref{dkdk}).

\n
By taking together the graded pieces of
the arrow $e: {\Bbb V} \to  {\Bbb V} [2](1)$, i.e. set $e':= \sum e_p$,
with $e_p: V_p \to V_{p+2}(1)$,
we obtain  
$e' \in L^{[0], \{2\}}_{(1)} (\tilde{\Bbb V}_*,\tilde{\Bbb V}_*)$ which, in turn,
induces the  homogeneous degree two arrow:
\beq\la{dfgt}
e \, : \, 
L^{[n]}_{(m)} ( \tilde{\Bbb V}_*,\tilde{\Bbb V}_* ) 
\lorw L^{[n]}_{(m+1)} ( \tilde{\Bbb V}_*,\tilde{\Bbb V}_* ),
\qquad 
u \longmapsto e(u):= e' \circ u - u \circ e'.
\eeq
There is a canonical arrow  of homogeneous degree $-2$ (this is where we need denominators (\ci{delseattle}, p.121):
\beq\la{dfgpp}
f \, : \, 
L^{[n]}_{(m)} ( \tilde{\Bbb V}_*, \tilde{\Bbb V}_* ) \lorw 
L^{[n]}_{(m-1)} ( \tilde{\Bbb V}_*,\tilde{\Bbb V}_* ).
\eeq
The arrows $(h,e,f)$ in (\ref{efbgt}), (\ref{dfgt}) and (\ref{dfgpp}) form
an ${\rm sl}_2 (\rat)$-triple turning the
 rational vector spaces
\beq\la{trfe}
L^{[n]}_{j} \, : \; \bigoplus_{d \in \zed^{\rm even/odd}} L^{[n], \{d\}}_{j+d/2}
 ( \tilde{\Bbb V}_*,\tilde{\Bbb V}_*)
\eeq
into ${\rm sl}_2 (\rat)$-modules; in what above, $j$ is a fixed integer multiple
of $1/2$ and the sum is over the integers $d$ with fixed parity,
even if $j$ is integral, odd if $j$ is an half-integer. Recalling that
the sum is finite, for $|d| \leq 2r$, we have that the corresponding
HL statememt reads as follows: ($e^k$ the $k$-th iteration of $e$ (\ref{dfgt}))
\beq\la{rftyq}
\xymatrix{
e^k \, : \;  L^{[n], (-k)}_{ (j-k/2) } \ar[rr]^{\hskip 0.7cm\cong} &&
L^{[n], (k)}_{ (j+k/2) }
}.
\eeq

Let  
$\varphi\,: {\Bbb V}_* \cong {\Bbb V}$ 
be any good  splitting (\ref{vsplits})
 and let $\tilde{e}(\varphi)$
be 
the associated matrix of $e: {\Bbb V} \to {\Bbb V}[2](1)$ (\ref{rvpqcr}).
The degree $d$ homogeneous  part of $\tilde{e}(\varphi)$ satisfies:
\beq\la{wxcer}
\tilde{e}(\varphi) ^{\{d\}}: = \sum_{q-p =d} \tilde{e}(\varphi )_{pq} \, \in \, L^{[2-d], \{d\} }_{(1)}
\eeq 
and is subject to (\ref{rrttnnyu}) and (\ref{ebbuono}):  it is zero for every $d>2$ and
it is the obvious arrow  for $d=2$.

The {\em the third Deligne isomorphism
 associated with $({\Bbb V}, e)$}:
\beq\la{rfml}\xymatrix{
\phi_{\rm III}(e) \,: \,= \phi_{\rm III} \, : \, {\Bbb V}_* \ar[rr]^{\hskip 1.4cm \cong} 
&&
 {\Bbb V},
}
\eeq  
is the unique  good splitting subject to
the following conditions: ($e^{1-d}$ is the $(1-d)$-iteration of (\ref{dfgt})):
\beq\la{condtr}
 e^{ 1-d } \left( \tilde{e}(\varphi)^{ \{ d \} } \right) \,=\,0, \qquad 
 \forall d\leq 1.
 \eeq
 
  \bigskip
 Let us illustrate how,  any  good splitting 
 $\varphi : {\Bbb V}_* \cong {\Bbb V}$ can be modified
 recursively, via HL (\ref{rftyq}), to obtain
 a new good splitting  subject to (\ref{condtr}). 
 
 Let $d=1$.  The condition (\ref{condtr}) reads $\tilde{e}^{ \{1 \} }=0$. Set $\varphi_1: =  \varphi ({\rm id} + \psi^{\{-1\} })$,
 where $\psi^{ \{-1\} } \in L^{ [1], \{-1\}}_{(0)}$ is a variable arrow.
 We conjugate $e$ and obtain:
 \beq\la{sdleq}
 \left( { \rm id} + \psi^{ \{ -1\} } \right)^{-1} \circ e \circ  \left(
 {\rm id} + \psi^{ \{-1\}} \right) \equiv
  \tilde{e}(\varphi)^{ \{2\} } + e \left ( \psi^{ \{-1\}}\right)  \quad \mbox{modulo degree $\leq 0$.}
 \eeq
 Note that the last term on the left is in $L^{ [1], \{1 \} }_{(1)  }$.
 We take the degree $1$ part of the r.h.s of (\ref{sdleq}) and set it equal to 
 zero
 \beq\la{edvgt}
  e \left(   \psi^{ \{-1 \} } \right)  = -\tilde{e}(\varphi)^{ \{1 \}}
  \qquad 
  \left(\mbox{equality in $L^{ [1], \{1\} }_{ (1)}$}\right).
 \eeq
 The HL (\ref{rftyq}) ensures that such a $\psi^{ \{-1\}}$ exists and is unique.
 This determines $\varphi_1$.

  Let $d=0$.  The condition (\ref{condtr}) reads $e(\tilde{e}^{ \{0 \} })=0$. Set $\varphi_0: =  \varphi_1 ({\rm id} + \psi^{\{-2\} })$,
 where $\psi^{ \{-2\} } \in L^{ [1], \{-2\}}_{(0)}$ is a variable arrow.
 We conjugate $e$ and obtain
 \beq\la{sdleq0}
 \left( { \rm id} + \psi^{ \{ -2\} } \right)^{-1} \, e\, \left(
 {\rm id} + \psi^{ \{-2\}} \right) \equiv
  \tilde{e} (\varphi_1)^{ \{2 \} } + e \left ( \psi^{ \{-2\}}\right)  \quad \mbox{modulo degree $\leq -1$.}
 \eeq
 Note that the last term on the left is in $L^{ [1], \{0 \} }_{(1)  }$.
 We take the degree $0$ part of the r.h.s of (\ref{sdleq0}) and set it equal to 
 zero after application of $e$:
 \beq\la{edvgt0}
  e^2 \left(   \psi^{ \{-2 \} } \right)  = - e\left( \tilde{e}(\varphi_1)^{ \{0 \}} \right)
  \qquad 
  \left(\mbox{equality in $L^{ [1], \{2\} }_{ (2)}$}\right).
 \eeq
 The HL (\ref{rftyq}) ensures that such a $\psi^{ \{-2\}}$ exists and is unique.
 This determines $\varphi_0$.
 
 We repeat this procedure for all decreasing values of $d$ and,
 recalling that ${\Bbb V}$ has type $[-r,r]$, the procedure ends
 no later than $d= -2r$. 
 
 The unicity of the resulting arrow is verified easily
 as follows. Let $a,b$ be two good splittings subject to
 (\ref{condtr}).
 Set   
 \beq\la{ee446677}
 \xymatrix{
 c\, :=\,  b^{-1} a = {\rm Id} + \sum_{l \geq 1} c^{ \{ - l\} } \, : \;
{\Bbb V}_*  \ar[rr]^{\hskip 2.6cm \cong} && {\Bbb V}.
}
\eeq
We apply the procedure carried out  above to $b$, 
 modifying it to $b_1:= b ({\rm Id} + c^{ \{-1\}} )$. We have that $b_1\equiv
 a$ modulo degree $\leq -2$, so that, in view of the fact that $a$  
 also satisfies (\ref{condtr}), we must have that $c^{ \{-1\} }=0$. It follows that  
 $b \equiv a$, modulo degree $\leq -1$. We repeat this procedure and kill
 all the $c^{ \{ -l\}}$.

 \bigskip
In general,  $\phi_{\rm I} \neq \phi_{\rm III}$, however, by 
 \ci{delseattle}, Proposition 3.6 (easily adapted to the present context),
we have that: (see Lemma (\ref{delsplgr}) for the definition of $f_i$)
\beq\la{11002299}
\phi_{\rm III} (e)_{| {\Bbb P}_{-i} } = 
\phi_{\rm I} (e)_{| {\Bbb P}_{-i} } = f_i : {\Bbb P}_{i} \lorw  {\Bbb V}.
\eeq

\begin{rmk}\la{rfcde}
{\rm
Unlike the four previous splittings $\omega_{\rm I,II}(e)$ and $\phi_{\rm I, II}(e)$,
the construction leading to $\phi_{\rm III} (e)$ is self-dual
in the sense that  we have:
\beq\la{kdkdslfd}
 \left(  \left(  \phi_{\rm III} \left( e^o \right) \right)^o \right)^{-1} = \phi_{\rm III}
 (e).
\eeq
This is because an isomorphism $\varphi$  satisfies  condition (\ref{condtr}) 
if and only if $\varphi^o$ satisfies the analogous ``opposite" conditions.
}
\end{rmk} 

\begin{rmk}\la{ercvoq}
{\rm
In genereal, 
the five good splittings 
\beq\la{eeemr4455}
\omega_{\rm I} (e), \quad \omega_{\rm II} (e),
\quad 
\phi_{\rm I}(e), \quad \phi_{\rm II}(e), \quad
\phi_{\rm III}(e) 
\eeq
are
pairwise distinct; see Examples \ref{snzdiff} and \ref{rrrr}.
}
\end{rmk}

\begin{pr}\la{eunico}
If an $e$-good splitting {\rm (\ref{gspltz})}  exists,
then it is unique 
and it coincides with the  third Deligne isomorphism.
In this case we have:
\beq\la{ee3rr344}
\omega_I (e) = \omega_{\rm II} (e) = \phi_{\rm I}(e) =\phi_{\rm II}(e) = \phi_{\rm III}(e).
\eeq
\end{pr}
{\em Proof.} Let $\varphi$ be $e$-good. Then $\tilde{e}(\varphi)^{ \{ d \} }=0$
for every $d \leq 1$. It follows   that condition (\ref{condtr}) is met
 by $\varphi$, so that  $\varphi = \phi_{\rm III}(e)$ and we have proved uniqueness
 (see also Remark \ref{fgtrrt}).

\n
Assume that there  is an $e$-good splitting, which, by the above, must coincide with $\phi_{\rm III} (e)$. By Remark \ref{fgtrrt}, we have that $\phi_{\rm III} (e) =
\phi_{\rm I} (e)$.  Proposition \ref{eecc11} implies that $\phi_{\rm I} (e)=\phi_{\rm II}(e)$ (this equality can be also    seen  by using a duality argument similar to the one in (\ref{559911})).

\n
Let us  compare $\phi_{\rm I}(e)$ with $\omega_{\rm I}(e)$.
By Remark (\ref{rvbblmlm}), the two agree on ${\Bbb V}_{-r} \oplus {\Bbb V}_{r}$.

\n
By comparing the general description (\ref{vneor}) of $\ke\, \rho$
with the  formula (\ref{evmoq}) for the embedding $\phi_{\rm I}$, we
see, with the aid of (\ref{lalskdjf}), that
\beq\la{fgfghre}
\ke \, \rho = \left\{ v \in V\, | \,  v= \sum_{(i,j) \in I_r} e^j  f_{i} (p_{ij})
 \right\},
\eeq
where $I_r$ is the set of the indices  subject to $0 \leq j \leq i$ 
and to $(i,j)\neq (r,0), (r,r)$. 
It follows that $\ke \, \rho$ coincides with $\phi_{\rm I} (\sum_{I_r} P_{-i}(-j))$.
By projecting onto $\sum_{|p|\neq r} V_p$, we deduce that
$\phi_{\rm I} (e)$, restricted to $\sum_{|p|\neq r} V_p$ factors
through $\ke \, \rho \to V$. 
Now we repeat for $\ke \, \rho$, what we have done above  for $V$
and deduce, by descending induction on $r$, that 
\beq\la{ecqape}
\phi_{\rm I}(e) =\omega_{\rm I}(e).
\eeq
We conclude by using  a duality argument:
\beq\la{559911}
\omega_{\rm II} (e) = ((\omega_{\rm I} (e^o))^o)^{-1}=
((\phi_{\rm I} (e^o))^o)^{-1} = \phi_{\rm II} (e) = \omega_{\rm I}(e),
\eeq
where: the first equality is by definition; the second equality
follows by the fact that there is an $e$-good splitting for
$({\Bbb V}, e)$ if and only  there is an $e^o$-good splitting
for $({\Bbb V}^o, e^o)$ and we have proved that $\omega_{\rm I} = \phi_{\rm I}$;
the third equality is by definition; the final equality is (\ref{ecqape}).
This concludes the proof.  
\blacksquare

\begin{ex}\la{snzdiff}
{\rm
$V \cong \rat^3$
with basis $(v_{-2}, v_0, v_2)$, $e:
(v_{-2}, v_0, v_2) \longmapsto (v_0, v_2, v_2 )$,
%\beq\la{wscbt}
%V = \langle v_{-2}, v_0, v_2 \rangle,
%\eeq
\beq\la{exvfr}
  \qquad V_{\leq -2} =  V_{\leq -1} =\langle v_{-2} \rangle, \quad
   V_{\leq 0} =V_{\leq 1}  =\langle v_{-2}, v_0 \rangle, \quad \; V_{\leq 2} =V, 
 \eeq
%\[ e\, : \;
%(v_{-2}, v_0, v_2) \longmapsto (v_0, v_2, v_2 ), \qquad e=
%\left(
%\begin{array}{ccc} 
%0 & 0&0\\
%1 & 0 & 0 \\
%0& 1& 1
%\end{array}
% \right)
%\]
\beq\la{exvfrgr}
  \qquad V_{ -2}  =\langle [v_{-2}] \rangle, \quad
   V_{ 0}=\langle [v_0] \rangle, \quad \; V_{2} = \langle [v_2] \rangle.
 \eeq
The five splittings  associated with $e$ discussed in this paper are: 
\bit
\item
$\phi_{\rm I} (e) \, :  ( [ v_{-2}], [v_0], [v_{2}]) \longmapsto
(  v_{-2}, v_0, v_2)$;

\item
$\phi_{\rm II} (e) \, :  ( [ v_{-2}], [v_0], [v_{2}]) \longmapsto
(  v_{-2},  -v_{-2}+ v_0, -v_0+v_2)$;

\item
$\phi_{\rm III} (e) \, :  ( [ v_{-2}], [v_0], [v_{2}]) \longmapsto
(  v_{-2},  -\frac{1}{3}v_{-2}+ v_0, -\frac{2}{3}v_0+v_2)$;

\item
$\omega_{\rm I} (e) \, :  ( [ v_{-2}], [v_0], [v_{2}]) \longmapsto
(  v_{-2},  -v_{-2}+ v_0, v_2)$;

\item
$\omega_{\rm I} (e) \, :  ( [ v_{-2}], [v_0], [v_{2}]) \longmapsto
(  v_{-2},  -v_{-2}+ v_0, -v_{-2} +v_2)$.
\eit
A direct calculation, or Proposition \ref{eunico}, shows that there is no $e$-good splitting.
}
\end{ex}
\begin{ex}\la{rrrr}
{\rm  Here is a class of  examples from geometry where, unlike
the previous example,  $e$ is nilpotent.
Let $Y\times Z$ be the product of  nonsingular complex projective varieties
and 
let $r:= \dim_{\comp}Z$.
Let $V:= H(Y\times Z,\rat)= H(Y,\rat)\otimes H(Z, \rat) 
= \oplus_{r,s} H^r(Y,\rat) \otimes H^s(Z,\rat)$.
Let $F'_s V$ be the subspace spanned by the
elements of the form $y_rz_{\s}$ with $\s\leq s$.
Set $F:= F'[r]$; this way $(V,F)$ has type $[-r,r]$.
Let 
$\eta \in H^2(Y\times Z,\rat)$ be the first Chern class
of a line bundle on $Y\times Z$ which is ample when restricted to the fibers
of the projection onto $Y$. Denote by 
$e: (V,F) \to (V,F[2])$ the map $ v \mapsto \eta \cup v$.
By using the hard Lefschetz theorem
on $Z$, one sees directly that the HL  condition (\ref{wsx}) holds for $((V,F),e)$. 
In addition to the five splittings  considered in this paper,
we also have the K\"unneth splitting $\kappa$.  
In general, the six splittings are pairwise distinct. 
 The reader can verify this fact  directly by taking $Y=Z$ to be an elliptic curve
and the line bundle to be of the form $E\times {\zeta} + {\zeta} \times E +
{\frak P}$, where   $\zeta \in E$ is a point and ${\frak P}$ is a Poincar\'e bundle
(this is to ensure that $\tilde{e}(\kappa)^{ \{1\} } \neq 0$,
so that, according to (\ref{condtr}), we must have $\kappa \neq \phi_{\rm III} (e)$).
If we take $Y =\pn{1} \times \pn{2}$, we are lead to examples
where $\kappa =\phi_{\rm III} (e)$, but otherwise the 
isomorphisms of type $\phi$ and $\omega$ are pairwise distinct, and distinct from 
$\phi_{\rm III}(e)$.
If we take $Y=Z=\pn{1}$, then  we have $\kappa = \phi_{\rm III} (e)$,
$\phi_{\rm I}(e) =   \omega_{\rm I} (e)$ and $\phi_{\rm II} (e) = \omega_{\rm II}(e)$,
but we have  no further relation. In this case, if we take $\eta$ to be the class
of the  fiber
of the projection onto $Z$, then we have an $e$-good splitting.
In all the examples, for general $\eta$, there is no $e$-good splitting.
A highly non-trivial example, 
where there is a
good splitting is mentioned at the end of the introduction.
}
\end{ex}

\section{Appendix: a letter from P. Deligne}
P. Deligne has sent the author a letter commenting on an earlier draft of this paper.
The author is happy to include, with P. Deligne's kind permission, this letter in this appendix as it outlines the simple modifications necessary to obtain
the splittings of this note in a Tannakian (tensor product) context.

\bigskip
March 16, 2012
\smallskip

Dear de Cataldo,

\smallskip
Thank you!

\smallskip
Let ${\mathscr A}$ be an Abelian category with a shift functor $A \to A(1)$ as in your text.
Let ${\mathscr B}$ be the category of objects of $\mathscr A$ given with a finite increasing filtration $F$ and $e: A \to A[2](1)$ verifying HL.

\medskip
\n
{\bf Corollary.} $\mathscr B$ is an abelian category.

\medskip
\n
{\em Proof.} (a) for objects $(A,F,e)$ of $\mathscr B$ of type $[-r,r]$, with $r>0$,
the ``peeling off" 2.2
\[
A= F_r \oplus \ke\, \left( 
F_{r-1} \subseteq A \stackrel{e^r}\lorw A \lorw A/F_{r-1}
\right)
\oplus e^r \, (F_r)
\]
is functorial. By induction on $r$, it follows that

\smallskip
\n
(b) the splitting $\omega_{\rm I}$ is functorial.

\smallskip
\n
If $f: (A,F,e) \lorw (A', F', e')$ is a morphism, (b) implies that
${\rm Gr}^F \ke\, (f) \stackrel{\cong}\lorw \ke\, {\rm Gr}^F (f)$,
and dually for coKer. That $(\ke\, (f), F, e)$ is in $\mathscr B$ follows.
It is a kernel. Dually for cokernels.
Morphisms are strictly compatible with filtrations, hence ${\rm CoIm} (f)
\stackrel{\cong}\lorw \im\, (f)$.

\medskip
\n
{\bf Remark.} In your 2.2. you should assume $a<b$.

\medskip
From now on, all categories are assumed to be $\rat$-linear.

\medskip
Let ${\mathscr A'}, {\mathscr A''}$ and ${\mathscr A}$ be as above and let 
$\otimes : {\mathscr A'} \times {\mathscr A''} \lorw {\mathscr A}$
be an exact 
biadditive
functor, compatible with shifts: functorial isomorphisms
$A' (1) \otimes A'' \stackrel{\cong}\lorw (A'\otimes A'')(1)$,
and $A'\otimes A''(1) \stackrel{\cong}\lorw (A'\otimes A'') (1)$ are given,
and the two resulting functorial isomorphisms
$A'(1) \otimes A''(1) \lorw (A'\otimes A'') (2)$ coincide.
Let ${\mathscr B'},  {\mathscr B''}$ and ${\mathscr B}$ be the corresponding categories of 
triples $(A,F,e)$ as above. Given $(A',F',e')$ and $(A'', F'', e'')$, one defines
the filtration $F$ (resp. morphism $e$) for $A' \otimes A''$ by 
\[
F_p = \sum_{p'+p''=p} F'_p \otimes F''_p
\]
(so that ${\rm Gr}^{F'} (A') \otimes {\rm Gr}^{F''} (A'') \stackrel{\cong}\lorw {\rm Gr}^F (A)$), and 
\[
e = e' \otimes 1 + 1 \otimes e'\]
(so that the same formula holds for the graded $e,e', e''$).

\medskip
\n
{\bf Proposition.} {\em This $\otimes$ sends ${\mathscr B}', {\mathscr B}''$ to
${\mathscr B}$.}

\medskip
One has to check that for graded objects and morphisms of degree $2$, $\otimes$ preserves
HL.
In the graded case, HL for $(A^*,e)$ is equivalent to the existence of $f^*: A^* \lorw
A^*(-1)$ of degree $-2$ such that $[e,f]$ id multiplicaiton by $n$ in degree $n$. Stability
of this property is proved in the same way that tensor product of representations of Lie algebras are defined.

Suppose now that ${\mathscr A}$ is a Tannakian category and that the twist is the tensor
product
with an object $\rat(1)$ of rank one. The compatibility between $\otimes$
and shift we required amounts tothe symmetry automorphism of $\rat(1) \otimes
\rat(1)$ being the identity.

\medskip
\n
{\bf Corollary.} {\em If ${\mathscr A}$ is Tannakian, so is ${\mathscr B}$.}

\medskip
\n
{\em Proof.} If $\omega$ is a fiber functor on ${\mathscr A}$, then $(A,F,e) \mapsto
\omega (A)$ is a fiber functor on ${\mathscr B}$.

\medskip
In terms of actions on $SL(2)$, rather in terms of grading and of $e$ and $f$ above,
I prefer to state the characterisitc property of the (good) splitting $\phi_{\rm III}$
as follows: it is the filtered isomorphism, with graded the identity:
\[
u\, :\, {\rm Gr}^F (A) \lorw A
\]
such that, with $e: {\rm Gr}^F (A) \lorw {\rm Gr}^F (A)(1)$, and $f$ as above, if $X$ is defiend
by
\[
u^{-1} e u = e+ X,
\]
one has $[f,X]=0$. This characterization makes it clear that this splitting
is compatible with tensor products (in the sense of the proposition). It gives an equivalence of the category ${\mathscr B}$ of triples $(A,F, e)$ with the category of graded objects $A^*$, $0$ outside  of finitely many degrees, given with
\[
e\, : \, A^* \lorw A^*(1) \;\; \mbox{and}\;\; f\, :\, A^* \lorw A(-1)
\]
of degree $2$, resp. $-2$, with $[e,f]=n$ in degree $n$, and given with $X: (\oplus A^n)
\lorw (\oplus A^n) (1)$ such that $[f,X] =0$.

In the Tannakian case, this equivalence is compatible with $\otimes$
[for $X$'s, $\otimes$ is defined by $X= X'\otimes 1 + 1 \otimes X''$].

\medskip
Best,

\smallskip
P. Deligne.


\begin{thebibliography}{99}



\bibitem{htam} M.A. de Cataldo, L. Migliorini,  The Hodge Theory of Algebraic maps,
Ann. Scient. \'Ecole
Norm. Sup., 4e serie, t. 38, (2005), 693-750. 

\bibitem{htadt} M. de Cataldo, L. Migliorini,  Hodge-theoretic aspects of the decomposition theorem, {\em Algebraic Geometry}, Seattle 2005, Proceedings of Symposia in Pure Mathematics, Vol. 80.2, 2009.

\bibitem{decmigpf} M.A. de Cataldo, L. Migliorini, ``The perverse filtration and the Lefschetz Hyperplane Theorem," 
            Annals of Mathematics, Vol. 171, No. 3, 2010, 2089-2113. 

\bibitem{decpf2} M.A. de Cataldo, ``The perverse filtration and the Lefschetz Hyperplane Theorem, II,"J. Algebraic Geometry 21 (2012) 305-345.

\bibitem{absho} M.A. de Cataldo, L. Migliorini, 
``The projectors of the decomposition theorem
are aboslute Hodge," submitted.

\bibitem{bams} M.A. de Cataldo, L. Migliorini, ``The Decomposition Theorem and the topology of algebraic maps," 
            Bullettin of the A.M.S., Vol. 46, n.4, (2009), 535-633. 

\bibitem{dechaumig} M.A. de Cataldo, T. Hausel, L. Migliorini, ``Topology of Hitchin systems and Hodge theory of character varieties: the case of $A_1$,''
            Annals of Mathematics 175 (2012), 1329-1407.
            
\bibitem{delseattle} P. Deligne,  Decompositions dans la categorie d\'eriv\'ee, {\em Motives}(Seattle, WA, 1991), 115,
128, Proc. Sympos. Pure Math. 55, Part 1, Amer. Math. Soc., Providence, RI, 1994.


\end{thebibliography}
\end{document}